\documentclass[letterpaper, 10pt, conference]{ieeeconf}
\IEEEoverridecommandlockouts
 
\usepackage{amsmath, amssymb, amsfonts, mathtools}
\usepackage{dsfont}
\usepackage{algorithmic}
\usepackage[short]{optidef}
\usepackage{graphicx}
\usepackage{textcomp}
\usepackage{xcolor}
\usepackage{booktabs}
\usepackage{multirow}
\usepackage{pifont} %
\usepackage{comment}
\usepackage[hidelinks]{hyperref}
\usepackage{flushend}

\usepackage{acronym}

\newacro{OCP}{optimal control problem}
\newacro{MPC}{model predictive control}
\newacro{RTI}{real-time iteration}
\newacro{RL}{reinforcement learning}
\newacro{QP}{quadratic programming}
\newacro{SQP}{sequential quadratic programming}
\newacro{CLC}{closed-loop costing}
\newacro{MDP}{Markov decision process}
\newacro{PPO}{proximal policy optimization}
\newacro{NLP}{nonlinear programming}
\newacro{GGN}{generalized Gauss-Newton}
\newacro{KKT}{Karush-Kuhn-Tucker}
\newacro{ODE}{ordinary differential equation}
\newacro{PEPT}{policy-enhanced partial tightening}
\newacro{PT}{partial tightening}

\newtheorem{assumption}{Assumption}
\newtheorem{remark}{Remark}
\newtheorem{lemma}{Lemma}
\newtheorem{proof}{Proof}
\newtheorem{definition}{Definition}
\newtheorem{theorem}{Theorem}
\newtheorem{corollary}{Corollary}

\DeclareMathSymbol{\shortminus}{\mathbin}{AMSa}{"39}

\newcommand{\R}{\mathbb{R}}
\newcommand{\N}{\mathbb{N}}
\newcommand{\Z}{\mathbb{Z}}
\newcommand{\norm}[1]{\left\lVert#1\right\rVert}
\newcommand{\set}[1]{\mathcal{#1}}

\newcommand{\epsu}{\varepsilon_u}

\newcommand{\V}{V}
\newcommand{\Vf}{V_f}

\newcommand{\state}{s}

\newcommand{\statevec}{\mathbf{x}}
\newcommand{\controlvec}{\mathbf{u}}
\newcommand{\slackvec}{\boldsymbol{\sigma}}

\newcommand{\piMPC}{\pi_{\rm MPC}}
\newcommand{\piMPCtau}{\pi_{\mathrm{MPC}, \tau}}
\newcommand{\piRL}{\pi_{\theta}}

\def\BibTeX{{\rm B\kern-.05em{\sc i\kern-.025em b}\kern-.08em
T\kern-.1667em\lower.7ex\hbox{E}\kern-.125emX}}

\begin{document}

\title{\bf Rollout Then Optimize: A One-Step Newton Refinement of Learned Policies for Nonlinear Model Predictive Control}
\author{{Andrea Ghezzi$^{1}$, Rudolf Reiter$^{2}$, Katrin Baumg\"artner$^{1}$, Alberto Bemporad${^3}$,  Moritz Diehl$^{1,4}$}
\thanks{$^1$ Department of Microsystems Engineering (IMTEK), University of Freiburg,~79110 Freiburg,~Germany}
\thanks{$^2$ Robotics and Perception Group, University of Zurich, Switzerland}
\thanks{$^3$ IMT School of Advanced Studies, Lucca, Italy}
\thanks{$^4$ Department of Mathematics, University of Freiburg}
\thanks{Correspondent: \tt andrea.ghezzi@imtek.uni-freiburg.de}
\thanks{This research was supported by DFG via projects 504452366 (SPP 2364) and 525018088, by BMWK via 03EI4057A and 03EN3054B, and by the EU via ELO-X 953348. This work was also supported by the European Research Council (ERC), Advanced Research Grant COMPACT (Grant Agreement No. 101141351).}}
\maketitle

\begin{abstract}
We propose a computationally efficient rollout-then-optimize method to improve a learned control policy at deployment time.
A learned policy provides a nominal trajectory, which is refined online by a single Newton step implemented via a Riccati recursion within a model predictive control (MPC) scheme.
This refinement combines model knowledge with the learned policy at minimal additional computational cost.
We establish bounds on the approximation error of the learned policy relative to the MPC policy and show that one Newton step reduces the suboptimality of the learned rollout quadratically in the policy approximation error.
The proposed controller is validated in simulation on a constrained trajectory-tracking task for a quadcopter with nonlinear dynamics.
Results highlight that the Newton step significantly improves the learned policy, achieving performance close to a fully converged MPC solution while requiring roughly half of the computational time.
The code is available at \footnotesize{\url{https://github.com/aghezz1/rl-riccati}}.
\end{abstract}

\section{Introduction}\label{sec:problem_statement}

\Ac{MPC} and \ac{RL} are two prominent approaches for approximately solving infinite-horizon optimal control problems~\cite{Reiter2026}.
Both provide feedback policies that minimize long-term cost while respecting system dynamics and constraints.
However, they differ significantly in how this objective is achieved and in their computational trade-offs.
MPC repeatedly solves finite-horizon optimal control problems online and can explicitly enforce constraints.
To achieve high performance, long prediction horizons are often required, which leads to substantial computational effort at deployment.
In contrast, RL shifts most of the computational burden offline: once training is complete, the learned policy can be evaluated very efficiently.
However, obtaining high-performance and constraint-satisfying policies typically requires extensive training and large amounts of data, especially for nonlinear systems.
The motivation of this work is to bridge this gap.
We aim to retain the fast online evaluation of learned policies while recovering part of the performance and constraint-handling capabilities of MPC without solving a nonlinear program (\acsu{NLP}) at each control step.

We consider the control of nonlinear dynamical systems $s_{t+1} = f(s_t, u_t)$, with state $s_t \in \set{S}\subseteq\R^{n_s}$, control $u_t\in \set{U}\subseteq\R^{n_u}, f:\set{S} \times \set{U} \to \set{S}$, and $t\in\mathbb{N}$.
The optimal feedback law $\pi^*:\set{S} \to \set{U}$ is defined as the solution of the infinite-horizon \acf{OCP}
\begin{mini}
	{\pi}{\sum_{t=0}^{\infty} \ell(s_t,u_t)}{\label{ocp: infinite horizon}}{J^*(s)\coloneqq}
	\addConstraint{s_0 }{= s}
    \addConstraint{s_{t+1}}{= f(s_t,u_t),}{\; t \in \N}
    \addConstraint{u_{t}}{= \pi(s_t),}{\; t \in \N}
	\addConstraint{g(s_t, u_t)  }{\leq 0,}{\; t \in \N,}
\end{mini}
where the stage cost~$\ell: \set{S}\times\set{U} \to \R_{\geq 0}$ and constraints $g: \set{S} \times \set{U} \to \R^{n_g}$ are continuous functions.

\ac{MPC} approximately solves \eqref{ocp: infinite horizon} by minimizing the cost over a finite-horizon under an open-loop policy~\cite{Rawlings2017}.
In every control step $t \in \mathbb{N}$, \ac{MPC} solves the NLP
\begin{mini!}
    {\statevec, \controlvec}{\Vf(x_N) + \sum_{k=0}^{N-1} \ell(x_k, u_k)}{\label{ocp: basic}}{\V^0_N(s_t) \coloneqq \label{cf: basic}}
    \addConstraint{x_0}{=s_t\label{cns:equality1}}
    \addConstraint{x_{k+1}}{=f(x_k, u_k),}{\; k \in \Z_0^{N-1}\label{cns:equality2}}
	\addConstraint{g(x_k, u_k)}{\leq 0,}{\; k \in \Z_0^{N-1},\label{cns:inequality}}
\end{mini!}
with $\Z_a^b \coloneqq \{a, a+1, \dots, b\}$ and decision variables $\statevec \coloneqq (x_0, \dots, x_N)$, $\controlvec \coloneqq (u_0, \ldots, u_{N-1} )$.
The optimal solution of \eqref{ocp: basic} is denoted by $\statevec^*(\state_t), \controlvec^*(\state_t)$, the MPC policy is denoted by $\piMPC(\state_t)$ and corresponds to the value of the first optimal control $u_0^*(s_t)$.

Conversely, \ac{RL} seeks to directly optimize a parameterized policy $\piRL$, typically represented by a neural network~\cite{Bertsekas1996a}.
This is done through iterative updates based on observed costs and transitions, making it specifically suited for problems where myriads of transition samples can be obtained, e.g., via simulation.
Both value-based and policy-gradient methods can be employed, as well as actor-critic architectures such as PPO~\cite{schulman2017ppo}.

Several approaches aim at reducing the online computational burden of MPC while retaining high performance.
One direction consists in learning a policy offline, either via reinforcement learning or imitation of an expert controller, and deploying it directly at runtime~\cite{hertneck2018learning,Ghezzi2023b}.
While this enables fast feedback evaluation, constraint satisfaction and performance guarantees are typically difficult to obtain~\cite{Gu2024}.
A different option for fast yet well-performing \ac{MPC} is to reduce the horizon while improving the terminal cost function \(\Vf\).
Ideally, \(\Vf\) should reflect the optimal cost-to-go beyond the finite horizon, reducing the impact of the horizon truncation on policy performance.
This approximation becomes increasingly important as the horizon length decreases.
Approaches include offline RL to learn $\Vf$~\cite{Reiter2024d}, using \ac{RL} to directly tune a parameterized and structured \(\Vf\)~\cite{Seel2022}, using supervised learning to construct a convex quadratic \(\Vf\) based on optimal state-control trajectories obtained from the solution of long horizon \ac{MPC}~\cite{Abdufattokhov2024}.
When a suboptimal policy is available, its rollout can be used to approximate the value function, an idea known as closed-loop costing (CLC)~\cite{Nicolao1998} or rollout of a base policy~\cite[\S 6.1.3]{Bertsekas1996a}.

\paragraph*{Contribution}
In contrast to approaches that directly deploy a learned controller or approximate the terminal value function offline, we follow a rollout-then-optimize principle.
We use a learned policy to generate a state-control trajectory which serves as linearization point for constructing a convex quadratic program (\acsu{QP}).
A single Newton step, implemented as a Riccati recursion within a \ac{RTI} MPC scheme, is performed on this \ac{QP}.
A backtracking line-search ensures constraint satisfaction.
Under standard assumptions, we show that the proposed controller is stabilizing and achieves lower closed-loop cost than the learned policy alone.
We demonstrate its effectiveness on a constrained trajectory-tracking task for a quadrotor.

\section{MPC Refinement of a Learned Policy Rollout}\label{sec:algorithm}
The proposed control scheme consists of two phases: an \emph{offline} and an \emph{online} phase.

\paragraph*{Offline phase}
In this phase, we train the policy~$\piRL$ in simulation.
This requires access to a simulator of the system dynamics~$f$ and to a cost function~$\hat{\ell}$ that reflects the objective of \eqref{ocp: basic}, possibly augmented to encode the system constraints.

\paragraph*{Online phase -- preparation}
We adopt a \ac{SQP} approach to solve~\eqref{ocp: basic}~\cite{Nocedal2006}.
To achieve real-time performance, we rely on the \ac{RTI} scheme~\cite{Diehl2005}, which performs a single SQP iteration per sampling instant.
Each iteration therefore requires one linearization of \eqref{ocp: basic} and the solution of one \ac{QP}.
We introduce nonnegative slack variables $\sigma_k\in\R_{\ge0}^{n_g}$ to rewrite inequality constraints in~\eqref{cns:inequality} as equalities.
The QP solved in each RTI step, linearized around a nominal trajectory $(\bar{\statevec},\bar{\controlvec})$, reads
\begin{mini!}
    {\statevec ,\controlvec ,\slackvec}
    {\sum_{k=0}^{N-1}\!
      \ell_k^{\mathrm{q}}(x_k,u_k)
      +\ell_N^{\mathrm{q}}(x_N)}
    {\label{ocp: qp}}{}
    \addConstraint{x_0\!-\!s_t}{= 0}
    \addConstraint{x_{k+1}\!-\!A_k x_k\!-\!B_k u_k\!-\!c_k}{= 0,}{\;k\in\Z_0^{N-1}}
    \addConstraint{d_k\!+\!G_k^x x_k\!+\!G_k^u u_k\!+\!\sigma_k}{= 0,}{\;k\in\Z_0^{N-1}\label{cns:slacked equality}}
    \addConstraint{\sigma_k}{\ge 0,}{\;k\in\Z_0^{N-1}.}
\end{mini!}
Here the quadratic approximations of the stage and terminal costs are
\begin{equation}\label{eq:qp cost components}
\begin{aligned}
\ell_k^{\mathrm{q}}(w_k)
  &\coloneqq \ell(\bar{w}_k)
   +\nabla\ell(\bar{w}_k)^{\!\top}\!\delta w_k
   +\tfrac{1}{2}\delta w_k^{\!\top}H_k\delta w_k,\\
\ell_N^{\mathrm{q}}(x_N)
  &\coloneqq \Vf(\bar{x}_N)
   +\nabla\Vf(\bar{x}_N)^{\!\top}\!\delta x_N
   +\tfrac{1}{2}\delta x_N^{\!\top} H_N \delta x_N,
\end{aligned}
\end{equation}
with $\delta w_k\coloneqq(\delta x_k, \delta u_k)\coloneqq(x_k-\bar{x}_k,u_k-\bar{u}_k)$.
The Hessian matrices are defined as (omitting the linearization point for brevity)
\begin{equation}\label{eq:exact Hessian of QP}
    H_N \coloneqq \nabla^2\Vf(\bar{x}_N), \; \begin{bmatrix}
      H_{xx} & H_{xu} \\
		H_{xu}^\top & H_{uu}
	\end{bmatrix}\coloneqq H_k, \; k \in \Z_0^{N-1}.
\end{equation}
In case the exact Hessian is used, $H_{xx} = \nabla^2_{xx}\mathcal{L}(\bar{w}_k, \bar{\lambda}_k, \bar{\mu}_k)$, $H_{xu} = \nabla^2_{xu}\mathcal{L}(\bar{w}_k, \bar{\lambda}_k, \bar{\mu}_k)$, $H_{uu} = \nabla^2_{uu}\mathcal{L}(\bar{w}_k, \bar{\lambda}_k, \bar{\mu}_k)$.
Function $\mathcal{L}$ is the Lagrangian of~\eqref{ocp: basic}, and $\bar{\lambda}_k \in \R^{n_x}, \bar{\mu}_k\in\R^{n_g}$ the multipliers associated with equality and inequality constraints, respectively.
Hence, $\bar{\boldsymbol{\lambda}}\coloneqq(\bar{\lambda}_0, \dots, \bar{\lambda}_N)$ and  $\bar{\boldsymbol{\mu}}\coloneqq(\bar{\mu}_0, \dots, \bar{\mu}_{N-1})$ are required and can be obtained from the previous closed-loop iteration (see Remark~\ref{remark:dependency on dual variables} below).
The linearized dynamics and constraints in~\eqref{ocp: qp} use the quantities
\begin{equation}\label{eq: linearized dynamics - gradients}
	\begin{split}
		& A_k \coloneqq \tfrac{\partial f}{\partial x} (\bar{x}_k, \bar{u}_k), \; B_k \coloneqq \tfrac{\partial f}{\partial u} (\bar{x}_k, \bar{u}_k), \\
		& G_k^x \coloneqq \tfrac{\partial g}{\partial x}(\bar{x}_k, \bar{u}_k), \; G_k^u \coloneqq \tfrac{\partial g}{\partial u}(\bar{x}_k, \bar{u}_k),
	\end{split}
\end{equation}
and
\begin{equation}\label{eq: linearized dynamics - zero order terms}
	\begin{split}
		& c_k \coloneqq f(\bar{x}_k, \bar{u}_k) - A_k \bar{x}_k - B_k \bar{u}_k, \\
		& d_k \coloneqq g(\bar{x}_k, \bar{u}_k) - G_k^x \bar{x}_k - G_k^u \bar{u}_k.
	\end{split}
\end{equation}

\paragraph*{Online phase -- feedback}
At each closed-loop step $t\in\N$, we perform a rollout $\bar{x}_{k+1} = f(\bar{x}_k, \piRL(\bar{x}_k)), \bar{x}_0 = s_t, k\in\Z_0^{N-1}$, using the learned RL policy $\piRL$.
We collect the state-control trajectory in the vectors $(\bar{\statevec}, \bar{\controlvec})$ and use it as linearization point for constructing $\eqref{ocp: qp}$.
The \ac{KKT} system of \eqref{ocp: qp} has to be solved.
For a stage $k\in\Z_1^{N-1}$, the system reads
\begin{subequations}\label{eq:kkt}
\begin{align}
H_{xx}\delta x_k + H_{xu}\delta u_k
 +A_k^\top \lambda_{k+1}
 -\lambda_k + {G_k^x}^{\!\top}\mu_k &= 0, \label{eq:kkt_a}\\
H_{uu}\delta u_k + H_{ux}\delta x_k
 +B_k^\top \lambda_{k+1}
 + {G_k^u}^{\!\top}\mu_k &= 0, \label{eq:kkt_b}\\
A_k\delta x_k + B_k\delta u_k - \delta x_{k+1} &= 0, \label{eq:kkt_c}\\
G_k^x\delta x_k + G_k^u\delta u_k + \sigma_k &= 0, \label{eq:kkt_d}\\
\sigma_k\ge0,\;\mu_k\ge0,\;
\operatorname{diag}(\sigma_k)\mu_k &=\tau\mathbf{1}, \label{eq:kkt_e}
\end{align}
\end{subequations}
where $\tau>0$ and $\mathbf{1}$ is a vector of ones with appropriate dimension.
We directly consider the relaxed complementarity condition \eqref{eq:kkt_e}, as we employ an interior-point solver that avoids us to explicitly handle active-set changes.
Optimality of \eqref{ocp: qp} is obtained for $\tau \to 0$.
Equation \eqref{eq:kkt_e} renders \eqref{eq:kkt} a nonlinear system and may require multiple Newton steps for convergence.
To limit online computation, we linearize \eqref{eq:kkt_e}
\begin{equation}\label{eq:kkt5 linearized}
    \operatorname{diag}(\sigma_k) \bar{\mu}_k + \operatorname{diag}(\bar{\sigma}_k) \mu_k = \tau \mathbf{1},
\end{equation}
where $\bar{\sigma}_k, \bar{\mu}_k$ denote slack and multiplier values from the previous closed-loop iteration.
Thus, the system composed of equations \eqref{eq:kkt_a}-\eqref{eq:kkt_d} and \eqref{eq:kkt5 linearized} with $k\in\Z_0^{N}$ is linear and can be solved with a single Newton step.
Finally, we need to perform a backtracking line-search to ensure $\sigma_k, \mu_k \geq 0, k\in\Z_0^{N-1}$.

\begin{remark}[Dependency on dual variables]\label{remark:dependency on dual variables}
If the QP \eqref{ocp: qp} is solved with exact Hessian, one would need access to the dual variables $\bar{\boldsymbol{\lambda}}, \bar{\boldsymbol{\mu}}$ to construct the KKT system \eqref{eq:kkt}.
However, the rollout with the RL policy provides only the primal variables (state–control trajectory) but not the multipliers.
Since \acp{OCP} typically employ a least-squares stage cost, it is common to adopt the Gauss-Newton Hessian approximation which does not require explicit dual information.
A second issue is related to the linearized complementarity \eqref{eq:kkt5 linearized} which requires $\bar{\mu}_k, \bar{\sigma}_k, k\in\Z_0^{N-1}$.
An automatic way to obtain these quantities is to solve the first $\text{QP}(s_0)$ to optimality, then its primal-dual solution can be used in \eqref{eq:kkt5 linearized} for the consecutive closed-loop iterations.
Hence, for $s_t, t\ge 1$, we can approximately solve the QP \eqref{ocp: qp} by performing a single Newton step on equations \eqref{eq:kkt_a}-\eqref{eq:kkt_d} and $\eqref{eq:kkt5 linearized}$ followed by a backtracking line-search.
\end{remark}

\begin{remark}[Block structure and Riccati recursion]
The linearized KKT equations form a symmetric block-tridiagonal system due to the sequential dynamics~\eqref{eq:kkt_c}.
This structure enables a forward-backward (Riccati-type) factorization~\cite{Rao1998}, allowing high-performance solution of \acp{QP}, as implemented for instance in the
open-source solver \texttt{HPIPM}~\cite{Frison2020a}.
\end{remark}

\section{Theoretical Properties}\label{sec: theory}
We next formalize the theoretical justification of the proposed method.

\begin{assumption}[Continuity and compactness]\label{ass:A0}
	The functions $f: \set{S} \times \set{U} \to \set{S}$, $\ell: \set{S} \times \set{U} \to \R_{\geq 0}$ and $\Vf: \set{S} \to \R_{\geq 0}$ are continuous, $f(0, 0) = 0$, $\ell(0, 0) = 0$, $\Vf(0) = 0$.
	The set $\set{S}$ is closed and $\set{U}$ is compact.
    Moreover, $f$ and~$\ell$ are Lipschitz continuous on $\set{S}\times\set{U}$ with constants $L_f,L_\ell>0$, and~$\Vf$ is Lipschitz with constant~$L_V>0$.
\end{assumption}
Under Assumption~\ref{ass:A0} a solution of \eqref{ocp: basic} exists~\cite[Prop. 2.4]{Rawlings2017}.

\begin{assumption}[MPC policy is nominally robust]\label{ass:A1}
For every $s_0\in\set{S}$ and sufficiently small bounded disturbances $e_t,w_t$, the closed-loop system
\[
s_{t+1}
 = f\bigl(s_t,\piMPC(s_t)+e_t\bigr)+w_t,  \quad t\in\N,
\]
is nominally robustly asymptotically stabilizing (NRAS) the origin $(s,u)=(0,0)\in\set{S}$.
\end{assumption}
Sufficient conditions for Assumption~\ref{ass:A0} are given, e.g., in~\cite[Th. 2.19, Pr. 3.5]{Rawlings2017}, and rely on Lyapunov stability.
Throughout the analysis we consider a deterministic setting
($e_t=w_t=0$) and constant horizon~$N$.
As $N\!\to\!\infty$, the policy~$\piMPC^N$ converges pointwise to the infinite-horizon optimal policy~$\pi^*$~\cite{Rawlings2017}.

Now consider training a parameterized policy~$\piRL$ (e.g., by reinforcement learning).
For an infinite number of simulations $N_{\mathrm{s}}$, the policy~$\piRL$ converges pointwise to~$\pi^*$.
In practice, however, the number of training episodes, hyperparameter tuning, and early stopping are limited, so~$\piRL$ only approximates the MPC feedback family.
We formalize this approximation in the following assumption.

\begin{definition}[Family of shrinking-horizon MPC policies]
For a finite horizon with $N$ steps, the solution $(\statevec^*, \controlvec^*)$ of~\eqref{ocp: basic} corresponds to the closed-loop system
\begin{equation}
\label{eq:MPC_cl}
\state_{t+1}=f\bigl(\state_t,\piMPCtau^{N-t}(\state_t)\bigr), \quad t\in \Z_0^{N-1},
\end{equation}
where $\piMPCtau^{N-t}: \set{S}\!\to\!\set{U}$ denotes the optimal control law for the shrinking horizon MPC problem with $N-t$ remaining steps and with $\tau$ being the minimum value of the barrier parameter achieved by the interior-point method solving \eqref{eq:kkt}.
For brevity, we denote the policy $\piMPCtau^N$ as $\piMPC$.
\end{definition}

\begin{assumption}[Approximation of shrinking-horizon MPC]
\label{ass:A2}
Suppose there exist sets $\set{S}_N \supseteq \set{S}_{N-1} \supseteq \dots \supseteq \set{S}_0$ such that
the origin $(s, u) = (0, 0) \in \set{S}_0$ and the closed-loop system under the learned policy is $f(s, \piRL(s)) \in \set{S}_{i-1}, \forall s \in \set{S}_i$.
Moreover, the learned policy $\piRL:\set{S}\to\set{U}$ satisfies the uniform stage-wise error bound
\begin{equation}\label{eq:policy_err}
\|\piRL(s)-\piMPCtau^{i}(s)\|
   \le \varepsilon_u,
   \;\;
   \forall s\in\set{S}_i,
\end{equation}
for some $\varepsilon_u>0$.
In addition, the policy outputs satisfy the problem constraints, i.e., $\piRL(s)\in\set{U}$ and $g(s,\piRL(s))\le0$ for all $s\in\set{S}$.
\end{assumption}

Thanks to Assumption~\ref{ass:A2}, the behavior of the closed-loop system $s_{t+1}=f(s_t,\piRL(s_t))$ is directly linked to that of the nominal MPC closed-loop $s_{t+1}=f(s_t,\piMPC(s_t))$.

\begin{corollary}
Under Assumption~\ref{ass:A1}, the system $s_{t+1} = f\bigl(s_t,\pi_{\mathrm{MPC}}(s_t)+e_t\bigr)$ remains NRAS for all bounded perturbations $\|e_t\|\le\varepsilon_u$ and all $s_t\in\set{S}_N$.
Consequently, the closed-loop system $s_{t+1}=f(s_t,\piRL(s_t))$ is asymptotically stable for sufficiently small~$\varepsilon_u$.
\end{corollary}

\begin{remark}
    With some additional assumption it is possible to extend the stability results to tracking problems with time-varying references~\cite[\S 2.4.5]{Rawlings2017}.
\end{remark}

\begin{remark}[Interpretation of the approximation bound]
The scalar $\varepsilon_u$ in~\eqref{eq:policy_err} represents the worst-case deviation (supremum over all stages $t$ and states $s$) between  the learned and the MPC policies.
It therefore provides a conservative measure of the approximation quality.
In practice, $\varepsilon_u$ depends on the training accuracy of $\piRL$, the smoothness of $f$ and the cost functions, and on how much the control law $\piMPCtau^{N-t}$ vary with~$t$ for a given $\tau\ge0$.
If stage-wise error bounds $\varepsilon_u^{N-t}$ are available, the subsequent state- and cost-error recursions can be written with time-varying terms, yielding tighter (less conservative) estimates of the trajectory deviation and cost difference.
\end{remark}

We turn our attention to the numerical solution of \eqref{ocp: basic} and consider the following standard assumptions.
\begin{definition}
We define by $z$ the concatenation of primal and dual variables of $\eqref{ocp: basic}$.
Then, $z\coloneqq(\statevec, \controlvec, \boldsymbol{\lambda}, \boldsymbol{\mu})$, where $\boldsymbol{\lambda} \in \R^{Nn_x}$ and $\boldsymbol{\mu} \in \R^{(N-1)n_g}$ are the multipliers associated with equality constraints~\eqref{cns:equality1}-\eqref{cns:equality2} and inequality constraints~\eqref{cns:inequality}, respectively.
We denote by $z^*$ the optimal solution of \eqref{ocp: basic}.
Let $R_\tau(z)$ denote the KKT residual of a nonlinear interior-point~(IP) formulation, obtained by relaxing complementarity as $\operatorname{diag}(\mu)\,g(x,u)=\tau\mathbf 1$ with $\tau>0$.
\end{definition}

\begin{assumption}[MPC regularity for IP method]\label{ass:A3}
For all $s\in\set{S}_N$, the optimal solution $z^*$ of~\eqref{ocp: basic} satisfies the linear independence constraint qualification (LICQ) and the second-order sufficient conditions (SOSC).
For $\tau$ sufficiently small, $z^*(\tau)$ denotes the solution of the perturbed system $R_\tau(z)=0$, and the Jacobian $\nabla R_\tau(z^*(\tau))$ is nonsingular in a neighborhood of $z^*(\tau)$.
\end{assumption}

\begin{assumption}[Lipschitz KKT Jacobian]\label{ass:A4}
There exists $\omega>0$ such that for all $z_1,z_2$ in a neighborhood of $z^*(\tau)$,
\begin{equation}
\bigl\|\nabla R_\tau(z_1)-\nabla R_\tau(z_2)\bigr\| \le \omega\,\|z_1-z_2\|.
\end{equation}
\end{assumption}

\begin{lemma}[Quadratic contraction of exact Newton step]
\label{lem:newton}
If the initial guess $z^{0}$ satisfies $\|z^{0}-z^*(\tau)\|\le\rho$ for a sufficiently small radius~$\rho$, the exact Newton step applied to the perturbed residual
\begin{equation}
z^{1}=z^{0}
      -\bigl[\nabla R_\tau(z^{0})\bigr]^{-1}R_\tau(z^{0}),
\end{equation}
obeys
\begin{equation}\label{eq:newton_quad_tau}
\|z^{1}-z^*(\tau)\|
   \le \tfrac{\omega}{2}\,\|z^{0}-z^*(\tau)\|^2.
\end{equation}
\end{lemma}

\begin{proof}
Standard local-quadratic result for Newton's method on a smooth perturbed KKT system applies (see e.g.~\cite[Th. 3.5 and \S 19]{Nocedal2006}).
For brevity, we omitted the dependence of $z^0, z^1$ on the barrier parameter $\tau>0$ in the notation.
The contraction result holds uniformly for small~$\tau$.
\end{proof}

\begin{theorem}[Quadratic improvement]\label{thm:final}
Let $\V$ be the cost functional $\V(\statevec,\controlvec)=\Vf(x_N)+\sum_{k=0}^{N-1}\ell(x_k,u_k)$ which is Lipschitz with constant $L_J$ (consequence of Assumption \ref{ass:A0}).
Let $z^{1}$ be the KKT vector obtained after a single Newton step starting from $z^0$ computed with the RL rollout
\[
z^{1}=z^{0}-\bigl[\nabla R_{\tau}(z^{0})\bigr]^{-1}R_{\tau}(z^{0}),
\qquad \tau>0~\text{small}.
\]
Then
\begin{equation}
\bigl\|\V(\statevec^{1},\controlvec^{1})
      - \V(\statevec^{*},\controlvec^{*}) \bigr\|
   \le L_J\,C_N\,\varepsilon_u^{2},
\end{equation}
with $C_N>0$ a constant depending only on the problem data and
uniform for sufficiently small~$\tau$.
Hence, the suboptimality of the RL-based trajectory is reduced quadratically with respect to the uniform control approximation error
$\epsu$.
\end{theorem}

\begin{proof}

Let $\mathbf{s}^{\mathrm{MPC}}$ and $\mathbf{s}^{\mathrm{RL}}$ denote the trajectories with length $N$ generated by rollouts of $\piMPCtau^{N-t}, t\in\Z_0^{N-1}$ and $\pi_{\theta}$, respectively, starting from the same $s_0$.
Define $e_k=s^{\mathrm{RL}}_k-s^{\mathrm{MPC}}_k$, $e_0=0$.
By Lipschitz continuity of $f$ and the uniform policy bound \eqref{eq:policy_err},
\begin{equation}\label{eq:error_rec}
\|e_{k+1}\|\le L_f(\|e_k\|+\varepsilon_u),\qquad k\in\Z_0^{N-1},
\end{equation}
which yields the closed-form bound
\[
\|e_k\|\le G_k\varepsilon_u,
\quad
G_k=\begin{cases}
\dfrac{L_f^k-1}{L_f-1}, & L_f\ne1,\\[1mm]
k,&L_f=1.
\end{cases}
\]
Hence, the state-error after $N$ steps satisfies $\norm{\mathbf s^{\mathrm{RL}}-\mathbf s^{\mathrm{MPC}}} \le G_N\varepsilon_u$.
Because the dual variables depend Lipschitz-continuously on the primal ones (LICQ), there exists $c_d>0$ such that
\begin{equation}\label{eq:KKT_err}
\|z^{0}-z^*(\tau)\|\le c_d G_N\varepsilon_u ,
\end{equation}
where $z^0$ is the KKT vector corresponding to the RL rollout and $z^*(\tau)$ the optimal solution of~\eqref{ocp: basic}.
Invoking Lemma~\ref{lem:newton} (applied to the perturbed residual $R_{\tau}$), one Newton step yields, uniformly for small~$\tau$,
\begin{equation}\label{eq:newton_err}
\|z^{1}-z^*(\tau)\| \le \tfrac{\omega}{2}c_d^{2}G_N^{2}\varepsilon_u^{2} =:C_N\varepsilon_u^{2}.
\end{equation}
The uniformity in $\tau$ ensures that the cost bound holds independently of the particular barrier value used in the IP formulation.
Since the cost functional $V(\statevec,\controlvec)=V_f(x_N) +\sum_{k=0}^{N-1}\ell(x_k,u_k)$ is Lipschitz on the feasible set with constant $L_J$,
\[
\big\|V(\statevec^{1},\controlvec^{1})
      -V(\statevec^{*},\controlvec^{*})\big\|
   \le L_J\|z^{1}-z^{*}(\tau)\|
   \le L_JC_N\varepsilon_u^{2}.
\]
This shows that the suboptimality after one Newton correction is of order~$\mathcal O(\varepsilon_u^{2})$.

\end{proof}

This result formalizes that, starting from the rollout of the learned policy, one Newton-step for the linearized MPC improves the cost of the resulting trajectory quadratically with respect to the policy approximation error~$\varepsilon_u$.

\begin{remark}[Contraction of Gauss-Newton Hessian]
The Newton step in Lemma~\ref{lem:newton} assumes access to the exact Jacobian of the KKT residual, which requires the Lagrange multipliers.
Similarly to what highlighted in Remark \ref{remark:dependency on dual variables}, in MPC formulation we can usually adopt the Gauss-Newton Hessian approximation which does not require dual information.
When using this approximation, the resulting iteration no longer enjoys quadratic contraction: convergence becomes at best linear, with a convergence factor that depends on the magnitude of the neglected higher-order residual terms (see, e.g.,~\cite[Ch.~10]{Nocedal2006}).
Nevertheless, this Gauss-Newton step often provides a sufficiently accurate correction in practice.
\end{remark}

\section{Practical Extensions of the Basic Algorithm}\label{sec: extensions}
We present extensions of the basic algorithm proposed in Sec.~\ref{sec:algorithm} that are useful in practice.

\paragraph*{Soft constraints or barriers}
To avoid possible solver failures due to numerical issues, modeling errors, or unknown disturbances, inequality constraints \eqref{cns:inequality} can be relaxed and penalized properly in the cost.
Alternatively, inequalities \eqref{cns:inequality} can be formulated in the cost as logarithmic barriers as done in \cite{Zanelli2017a}.
However, barrier functions are not Lipschitz continuous and would therefore invalidate the Lipschitz assumptions used in the theoretical analysis above.

\paragraph*{Synthesizing terminal cost}
The Riccati recursion performed on the KKT system with linearized complementarity conditions of \eqref{ocp: qp} can be interpreted as a method for synthesizing approximate terminal costs.
Indeed, the by-product of the Riccati recursion is a positive-semidefinite Hessian matrix and a gradient of the value function $V_N^0(s_t)$.
The formulas of the Riccati recursion and corresponding value function are reported explicitly in~\cite{Zanelli2017a}.
This terminal cost interpretation of the Riccati recursion allows defining an MPC problem with two phases (or horizons) as in~\cite{Nicolao1998}.
The first phase, $k\in\Z_0^{M-1}, M<N$, called the \emph{control horizon}, corresponds to a standard MPC problem as in \eqref{ocp: basic}.
The second phase, $k\in\Z_{M}^N$, called the \emph{simulation horizon}, applies the proposed algorithm.
In the second phase, we generate a rollout using $\piRL$ starting from $\bar{x}_M$, construct a QP, and derive its KKT system with linearized complementarity.
Then, a Riccati recursion is performed backward from $N$ to $M$ to solve the KKT system.
The recursion yields a Hessian $H_M(\bar{x}_M)\succeq0$ and a gradient $h_M(\bar{x}_M)$ which compose the quadratic terminal cost term used in the first-phase MPC problem.

\begin{remark}
Both the soft-constrained and the two-phase formulations can be useful in practice, especially when the RL policy is still being learned and its performance is not yet satisfactory.
Extending the first (control) horizon and shrinking the second (simulation) horizon makes the overall MPC problem approach closer to a standard single-phase formulation.
Thanks to its limited online computation, the proposed approach is attractive for RL training frameworks employing MPC-in-the-loop, where thousands of MPC evaluations are typically needed~\cite{Reiter2026}.
\end{remark}

\paragraph*{Real-time implementation}\label{parag:rti}
In a RTI scheme it is customary to divide the computations into two phases, \emph{preparation} and \emph{feedback}~\cite{Diehl2005}.
During the preparation phase, the KKT matrix of the QP~\eqref{ocp: qp} can be factorized without requiring the latest state measurement.
Once the new state is available, the \emph{feedback} phase starts and the QP is solved to compute the control law.
In the proposed algorithm, before assembling the QP, we perform the $N$-step rollout with the learned policy~$\piRL$, which in principle should start from the current measured state~$s_t$.
However, performing this rollout only after obtaining~$s_t$ shifts all computations to the feedback phase and may introduce an undesirably large delay between the measurement of~$s_t$ and the application of the control~$\pi_{\mathrm{MPC}}(s_t)$.
A practical way to minimize feedback time is to compute the rollout in advance, starting from~$x_1^*(s_{t-1})$, i.e., the one-step predicted optimal state in the previous closed-loop iteration.

\section{Trajectory tracking with a quadcopter}\label{sec: simulations}
We test the proposed scheme considering the control of a nano-quadcopter that has to track a lemniscate trajectory in the 3D space.
The results are obtained via numerical simulations using the quadcopter environment of \small \texttt{safe-control-gym}\normalsize~\cite{Yuan2021safecontrolgym}.
The MPC problems described next are formulated and solved using \texttt{acados}~\cite{Verschueren2021} with \texttt{HPIPM} as QP solver~\cite{Frison2020a}.
All computations are performed on a MacBook with Intel(R) Core(TM) i7-8559U~CPU~@~2.70GHz and 16 GB of memory running macOS 15.1.

\subsection{Quadcopter modeling}
The quadcopter state is defined as $s \coloneqq (\mathbf{p}, \mathbf{v}, \boldsymbol{\Psi}, \boldsymbol{\omega})$, where  $\mathbf{p} \in \R^3$ is the position of the quadcopter center of mass (CoM) in the inertial frame, $\mathbf{v} \in \R^3$ is the linear velocity of the CoM in the world frame, $\boldsymbol{\Psi} \in \R^3$, is the orientation in the body frame expressed as roll, pitch, yaw, and $\boldsymbol{\omega} \in \R^3$ contains angular velocities in the body frame with respect to the inertial frame.
The control of the quadcopter corresponds to the thrust of single propellers, $u \coloneqq (\tau_1, \tau_2, \tau_3, \tau_4) \in \R^4$.
The dynamics of the drone are described by the \ac{ODE}
\begin{equation}\label{eq: drone ode}
	\dot{s} \!= \!f(s, u) \!=\!
	\left(\begin{matrix}
		\mathbf{v}  \\
		\frac{1}{m} \mathbf{R}_\mathrm{bw} \mathbf{e}_z \sum_{i=1}^{4}\tau_i - \mathbf{e}_z g \\
		\mathbf{T} \boldsymbol{\omega} \\
		J^{-1} \!\left(
			\mathbf{M}(\tau_1, \tau_2, \tau_3, \tau_4)
			- \boldsymbol{\omega} \times J \boldsymbol{\omega}
			\right)
	\end{matrix}\right),
\end{equation}
where $ \mathbf{R}_\mathrm{bw} \in \R^{3\times3}$ is the rotation matrix from the body frame to the world frame, $\mathbf{e}_z \in \R^3$ is the unit vector for the $z$-direction of the world frame, $g$ is the gravity acceleration, $\mathbf{T}\in \R^{3\times3}$ is a matrix that maps angular velocity of the body frame into the time derivative of the orientation, $J\in \R^{3\times3}$ is the inertia matrix of the quadcopter, and $\mathbf{M}\in \R^{3}$ is a vector containing the angular momentum generated by the propellers.
We adopt SI units.
For more details and parameter values of the model see~\cite{Yuan2021safecontrolgym,luis2016design}.

\subsection{MPC and RL policies}
The MPC policy to accomplish the tracking task solves in closed-loop iteration $t\in \N$ the following \ac{NLP}
\begin{mini}
	{\statevec, \controlvec}{\sum_{k=0}^{N-1}
		\delta z_k^\top H_k \delta z_k +
        \delta z_N^\top H_N \delta z_N
        }{\label{ocp: drone example}}{}
	\addConstraint{x_0 - s_t}{=0}
	\addConstraint{x_{k+1} - F(x_k, x_k)}{= 0,}{\;k \in \Z_{[0, N-1]}}
	\addConstraint{u_\mathrm{lb} \leq u_k}{\leq u_\mathrm{ub}, }{\;k \in \Z_{[0, N-1]}}
    \addConstraint{x_\mathrm{lb} \leq x_k}{\leq x_\mathrm{ub}, }{\;k \in \Z_{[0, N]},}
\end{mini}
where $N=40$, $\delta z_k \coloneqq (x_k - x_{t+k}^o, u_k - u_{t+k}^o)$, $\delta z_N \coloneqq (x_N - x_{t+N}^o)$, and $x_{t+k}^o, u_{t+k}^o$ are references to track.
The quadratic cost is defined by the matrix $H_N \coloneqq Q$, and by the block diagonal matrix $H_k \coloneqq \mathrm{diag}(Q, R)$, with $Q = \text{diag}(5, 5, 5, 0.1 \cdot \mathbf{1}_8), \quad R = \text{diag}(0.1 \cdot \mathbf{1}_4)$, and $\mathbf{1}_n$ is a vector of ones with dimension $n$.
Function $F: \R^{n_x}\times\R^{n_u} \to \R^{n_x}$ corresponds to the time-discretization of the \ac{ODE}~\eqref{eq: drone ode} with an explicit Runge-Kutta integrator of order 4, and sampling time $t_\mathrm{d}=0.02$~s.
The position reference consists of a 3D lemniscate with a period of $T = 5$~s and frequency $\rho = \frac{2\pi}{T}$, its derivative constitutes the reference for the velocity.
The lemniscate equations and numerical values for the bounds $u_\mathrm{lb}, x_\mathrm{lb}, u_\mathrm{ub}, x_\mathrm{ub}$ are available online\footnote[1]{\url{https://github.com/aghezz1/rl-riccati}}.
For the orientation and angular velocity of the quadcopter we consider the reference to be zero.
Note that the resulting reference might be dynamically infeasible for the quadcopter.
As a reference for the controls, we consider the thrust required for hovering, $u^o_{t+k} = 0.06615 \cdot \mathbf{1}_4 \text{ N }, k \in \Z_{[0, N-1]}$.
Finally, we adopt a Gauss-Newton Hessian approximation.

\paragraph*{Reinforcement-learning policy}
We adopt the RL policy obtained via PPO that is available in \small\texttt{safe-control-gym} \normalsize for the specified quadcopter model and tracking task.
The policy input consists of the current state $s_t$ and the state reference to track one-step forward, i.e., $x^o_{t+1}$.
The RL policy has been trained by minimizing the stage cost of the MPC~\eqref{ocp: drone example}.
A training episode of the RL policy is truncated in case of constraint violations and incurs a higher cost.
Hence, we promote safer but potentially suboptimal trajectories.
From the stochastic PPO policy parameterized as Gaussian, we consider the predicted mean value to obtain a deterministic policy, denoted here as PPO-RL.
The PPO-RL policy is a fully-connected neural network with one hidden layer with 128 elements and $\tanh$ activation functions.

\subsection{Comparison}
We compare four policies: \emph{RTI}, \emph{Riccati-RTI}, \emph{Riccati-RL}, and \emph{PPO-RL}.
These correspond respectively to solving~\eqref{ocp: drone example}
(i) with a standard RTI method,
(ii) with an RTI method limited to a single iteration,
(iii) with the proposed method described in Sec.~\ref{sec:algorithm},
and (iv) directly with the trained PPO policy.

For the RTI controller, the quadratic programs~\eqref{ocp: qp} are solved by iterating on the smoothed nonlinear KKT system~\eqref{eq:kkt}, starting from an initial barrier parameter~$\tau_t=1$ and driving it to $\tau_{\min}=10^{-8}$ at the optimum.
The Riccati-RTI and Riccati-RL controllers perform only one Newton step on the primal-dual system \eqref{eq:kkt} with the linearized complementarity condition~\eqref{eq:kkt5 linearized}.
Regarding the initial point of each Newton step, primal variables in Riccati-RTI are initialized using the previous solution, and for the very first step we use the tracking reference.
Instead, for Riccati-RL, the primal variables are initialized via a rollout of the PPO policy.
Dual variables are initialized using the previous solution for both approaches, the initial value of the barrier parameter is determined as $\tau_t = \| \text{diag}(\sigma_{t-1}) \mu_{t-1} \|_1$.

\begin{figure*}
    \centering
    \includegraphics{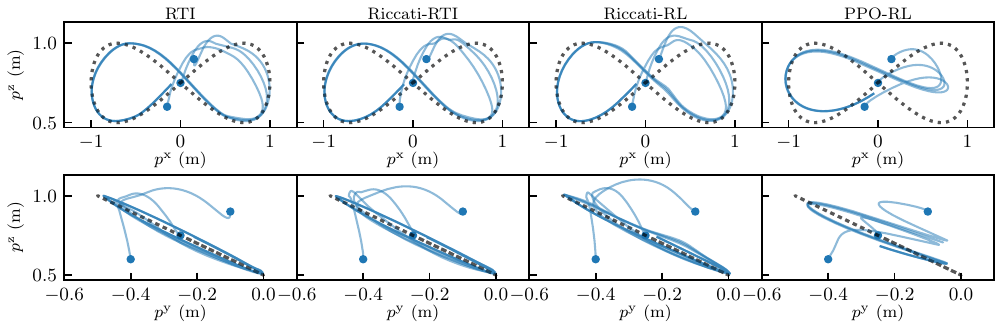}
    \includegraphics[trim={0, 0, 0, 0.355cm}, clip]{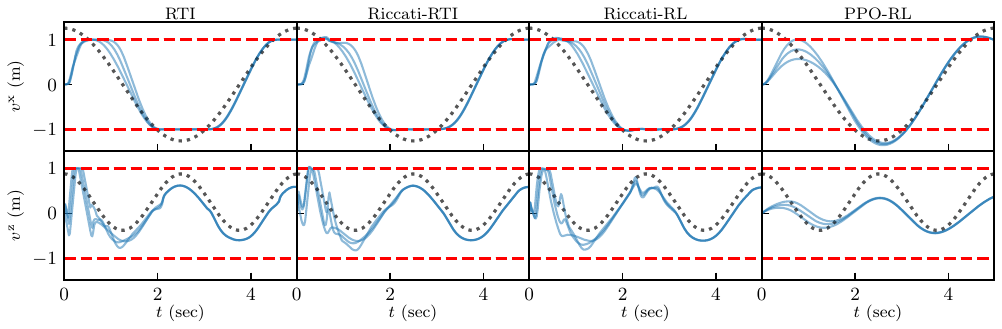}
    \caption{Comparison of trajectory tracking performance for four control strategies: RTI, Riccati-RTI, Riccati-RL (\emph{the proposed method}), and PPO-RL.
    The first two rows show the position trajectories in the space $p^\mathrm{x},p^\mathrm{z}$ and $p^\mathrm{y},p^\mathrm{z}$, respectively.
    The blue dots correspond to the initial position of each simulation.
    The bottom two rows present the time evolution of velocity $v^\mathrm{x}$ and $v^\mathrm{z}$, respectively.
    The black dashed lines indicate the reference trajectories, and the red dashed lines denote constraint limits.}
    \label{fig:state traj}
\end{figure*}

Three closed-loop simulations with duration $T$ are performed starting from a hovering condition and from three initial positions shifted along the lemniscate by $\pm 0.15$ m with respect to its center.
We compare tracking accuracy, constraint violations, and computational time.
Figure~\ref{fig:state traj} shows the time trajectories of positions and linear velocities along the $\mathrm{x}$- and $\mathrm{z}$-axes.

\begin{table}
    \vspace{0.15cm}
	\caption{Comparison of closed-loop cost}
    \vspace{-0.2cm}
	\centering
    \scriptsize
    \begin{tabular}{@{}lccccccc@{}}
    \toprule
    \multirow{2}{*}{Approach} &
    \multicolumn{3}{c}{Tracking Cost} &
    \multicolumn{3}{c}{Constraint Violation Cost} \\
    \cmidrule(lr){2-4} \cmidrule(lr){5-7}
     & Sim. 1 & Sim. 2 & Sim. 3 & Sim. 1 & Sim. 2 & Sim. 3 \\
    \midrule
    RTI                 & 0.086 & 0.073 & 0.145 & 0 & 0 & 0 \\
    Riccati-RTI         & 0.095 & 0.083 & 0.179 & 0.027 & 0.012 & 0.033 \\
    Riccati-RL          & 0.101 & 0.087 & 0.177 & 0.022 & 0.016 & 0.022 \\
    PPO-RL              & 0.181 & 0.177 & 0.229 & 0.155 & 0.144 & 0.159 \\
    \bottomrule
    \end{tabular}
    \label{tab:costs}
    \vspace{-0.4cm}
\end{table}

Table~\ref{tab:costs} reports the average closed-loop cost for each simulation.
The cost consists of a tracking term, identical to the quadratic stage cost in~\eqref{ocp: drone example}, and a penalty for constraint violations, defined as $\mathcal{J}_{\mathrm v} = \frac{1}{N_\mathrm s}\sum_{i=1}^{N_\mathrm s}\!\norm{\max(0,h(z))}_1$, where $N_\mathrm s = T/t_\mathrm d$ and $h(z)\le0$ collects the state and input bounds.
As expected, RTI achieves the lowest tracking cost in all scenarios, with no constraint violations, and thus serves as the performance reference.
Restricting RTI to a single iteration (Riccati-RTI) degrades the performance, yielding higher tracking cost and some constraint violation cost.
The purely learned PPO-RL controller incurs in about double the tracking cost compared to RTI, and it exhibits constraint violations, concentrated mainly on the linear velocity~$v^{\mathrm x}$.
This is why we omitted plots of other states or controls.
Finally, the proposed Riccati-RL scheme is able to recover almost the same performance of RTI, dramatically improving both tracking and violation cost compared to PPO-RL, confirming the benefit of the optimization-based refinement step.

\begin{table}
    \vspace{0.15cm}
	\caption{Average and maximum runtime for one control step}
    \vspace{-0.2cm}
	\centering
    \scriptsize
    \begin{tabular}{@{}lccccccccc@{}}
    \toprule
    \multirow{2}{*}{Approach} &
    \multicolumn{2}{c}{Initialization (ms)} &
    \multicolumn{2}{c}{Runtime (ms)} &
    \multicolumn{2}{c}{Feedback Time (ms)} \\
    \cmidrule(lr){2-3} \cmidrule(lr){4-5} \cmidrule(lr){6-7}
     & Avg. & Max. & Avg. & Max. & Avg. & Max. \\
    \midrule
    RTI            & -    & -      & 1.86 & 4.55 & 1.00 & 3.58 \\
    Riccati-RTI    & -    & -      & 1.00 & 1.66 & 0.17  & 0.47    \\
    Riccati-RL     & 11.56 & 15.09 & 1.11 & 2.27 & 0.19 & 0.31 \\
    PPO-RL         & -    & -      & 0.16 & 0.46 & -    & -    \\
    \bottomrule
    \end{tabular}
    \label{tab:runtime}
    \vspace{-0.4cm}
\end{table}

Table~\ref{tab:runtime} summarizes the computational performance in terms of initialization time (``Initialization''), total solver runtime (``Runtime''), and \emph{feedback} time (see Sec.~\ref{sec: extensions}a).
For the RTI-based schemes, ``Runtime'' is the sum of the \emph{preparation} and \emph{feedback} phases.
The ``Initialization'' time accounts for performing the RL rollout and setting the resulting state-control trajectory as the initial guess for \texttt{HPIPM}.
The \emph{preparation} phase is identical for all optimization-based approaches since they assemble the same QP~\eqref{ocp: qp}.
Both Riccati-RTI and Riccati-RL exhibit a much shorter \emph{feedback} time than full RTI because they perform only a single Riccati recursion.
As anticipated, PPO-RL achieves the smallest overall runtime because it evaluates only a neural network.
The measured initialization time of Riccati-RL is relatively large, mainly due to Python overhead in the implementation of the policy rollout.
In practice, the rollout involves the evaluation of a small fully-connected neural network and the integration of the system dynamics with a simple explicit Runge-Kutta scheme over~$N$ steps.
Therefore, its execution time is expected to decrease by at least one order of magnitude in compiled code.
Despite this overhead, all methods remain well below the 20\,ms sampling interval, guaranteeing realtime feasibility.

\section{Conclusions and Outlook}
This work demonstrated that learned control policies can be effectively improved online through a single Riccati-based Newton correction executed within a real-time-iteration MPC framework.
We established theoretical bounds on the approximation error between the learned and MPC policies and proved that a single correction step reduces the trajectory suboptimality quadratically with the learned-policy error magnitude.
Simulations on a quadcopter tracking problem confirmed that the Riccati-RL controller significantly improves the learned policy, achieving performance close to RTI at lower computational cost.
Future work will focus on a high-performance implementation of the policy rollout to minimize initialization overhead, and on experimental validation on real hardware to confirm the predicted computational advantages and closed-loop performance.

\bibliography{syscop,additional_lib}

\end{document}